\documentstyle[12pt,amssymb,graphicx]{article}

\begin{document}
\renewcommand{\thefootnote}{}
\date{}

\def\thebibliography#1{\noindent{\normalsize\bf References}
 \list{{\bf
 \arabic{enumi}}.}{\settowidth\labelwidth{[#1]}\leftmargin\labelwidth
 \advance\leftmargin\labelsep
 \usecounter{enumi}}
 \def\newblock{\hskip .11em plus .33em minus .07em}
 \sloppy\clubpenalty4000\widowpenalty4000
 \sfcode`\.=1000\relax}

\title{\vspace*{-1cm}
{%\normalsize
\large 
Classification of $n$-component Brunnian links up to 
$C_n$-move}
}

\author{
{\normalsize Haruko Aida MIYAZAWA
\thanks{The first author is partially supported by the 21st COE program 
\lq\lq Constitution of wide-angle mathematical
basis focused on knots''.}}
\\
{\small Research Institute for Mathematics
and Computer Science,}\\[-1mm]
{\small Tsuda College}\\[-1mm]
{\small Kodaira, Tokyo 187-8577, Japan}\\[-1mm]
{\small e-mail: aida@tsuda.ac.jp}\\[3mm]
{\normalsize Akira YASUHARA }\\
{\small Department of Mathematics, Tokyo Gakugei University}\\[-1mm]
{\small Nukuikita 4-1-1, Koganei, Tokyo 184-8501, Japan}\\[-1mm]
{\small e-mail: yasuhara@u-gakugei.ac.jp}\\[3mm]
%{\small {\bf Please send all correspondances to 
%the second author.}}
}

\maketitle

\vspace*{-5mm}  
\baselineskip=15pt
{\small 
\begin{quote}
\begin{center}A{\sc bstract}\end{center}
\hspace*{1em} We give a classification of $n$-component links 
up to $C_n$-move. In order to prove this classification, 
we characterize Brunnian links, and have that 
a Brunnian link is ambient isotopic to a band sum of 
trivial link and Milnor's links.
\end{quote}}

\footnote{
%; Secondary 57M27 %Invariants of knots and 3-manifolds
}
%\footnote{{\em Short Running Title}: }

\renewcommand{\thefootnote}{*}

\baselineskip=15pt

\noindent
{\bf 1. Introduction} 

\bigskip
A {\em tangle} $T$ is a disjoint union of properly embedded 
arcs in the unit $3$-ball $B^{3}$. 
A tangle $T$ is {\em trivial} if there 
exists a properly embedded disk in $B^3$ that contains $T$. 
A {\em local move} is a pair of trivial tangles 
$(T_{1},T_{2})$ with $\partial T_{1}=\partial T_{2}$ 
such that for each component $t$ of $T_1$ there exists 
a component $u$ of $T_2$ with $\partial t=\partial u$.
Two local moves $(T_{1},T_{2})$ and $(U_{1},U_{2})$
are {\em equivalent} 
if there is an orientation preserving 
self-homeomorphism $\psi :B^{3}\rightarrow B^{3}$ such that $\psi (T_{i})$ 
and $U_{i}$ are ambient isotopic in $B^3$ relative to $\partial
B^{3}$ for $i=1,2$. 
The definition of a local move follows from \cite{T-Y2}, \cite{T-Y}. 

A {\em $C_{1}$-move} is a local move $(T_1,T_2)$ as illustrated in Fig. 1. 
Suppose that a $C_k$-move is defined.  
Let $(T_{1},T_{2})$ be a 
$C_k$-move, $t_{1}$ a component of $T_1$ and $t_2$ a component of $T_2$ 
with $\partial t_{1}=\partial t_{2}$. Let $N_i\ (i=1,2)$ be regular 
neighbourhoods of $t_i$ in $(B^3-T_i)\cup t_i$ such
that $N_{1}\cap \partial  B^{3}=N_{2}\cap \partial B^{3}$. 
Let $\alpha$ be a disjoint union of properly 
embedded arcs in $B^{2}\times I$ as illustrated in Fig. 2.
Let $\psi_{i}:B^{2}\times I\rightarrow N_{i}$ be a homeomorphism 
with $\psi_{i}(B^{2}\times \{ 0,1\} )=N_{i}\cap \partial B^{3}$ for $i=1,2$. 
Suppose that $\psi_{1}(\partial \alpha )=\psi_{2}(\partial \alpha )$ and 
that $\psi_{1}(\alpha )$ and $\psi_{2}(\alpha )$ are ambient 
isotopic in $B^{3}$ 
relative to $\partial B^3$. Then we say that a local move 
$((T_{1}-t_{1})\cup \psi_{1}(\alpha ), (T_{2}-t_{2})\cup \psi_{2}(\alpha ))$ 
is a {\em $C_{k+1}$-move}. 
Note that, for each natural number $k$, 
there are finitely many $C_{k}$-moves up to equivalence. 
It is easy to see that if $(T_1,T_2)$ is a $C_k$-move, 
then $(T_2,T_1)$ is equivalent to a $C_{k}$-move.
The {\em one-branched $C_{k}$-move} is a local move as 
illustrated in Fig. 3. 
It is shown that a $C_k$-move is generated by 
the one-branched $C_k$-move \cite{Habiro2}, \cite{T-Y}. 
So the one-branched $C_k$-move is often called the $C_k$-move. 
The definition of $C_k$-move is due to Habiro \cite{Habiro2}, \cite{Habiro1}.

\begin{center} 
\begin{tabular}{cc}
\includegraphics[trim=0mm 0mm 0mm 0mm, width=.35\linewidth]
{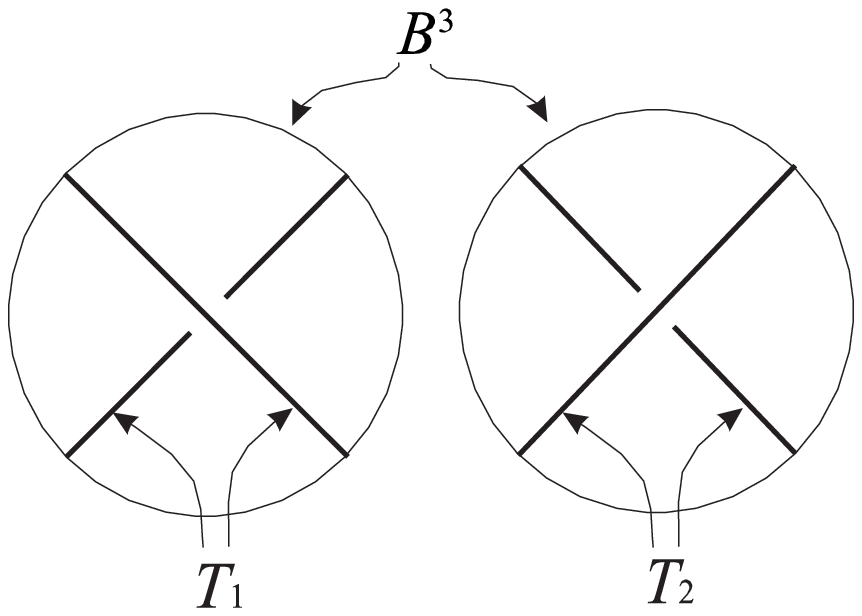} \hspace*{1in} &
\includegraphics[trim=0mm 0mm 0mm 0mm, width=.25\linewidth]
{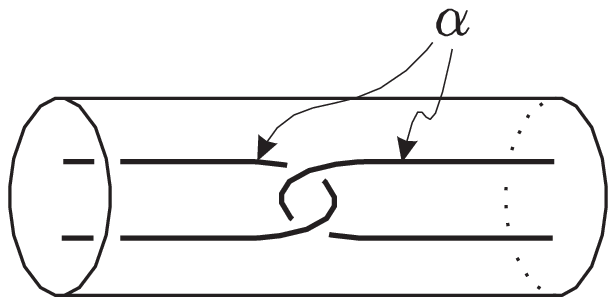}\\
Fig. 1 \hspace*{1in} &
Fig. 2
\end{tabular}
\end{center}

\medskip
\begin{center}
\includegraphics[trim=0mm 0mm 0mm 0mm, width=.6\linewidth]
{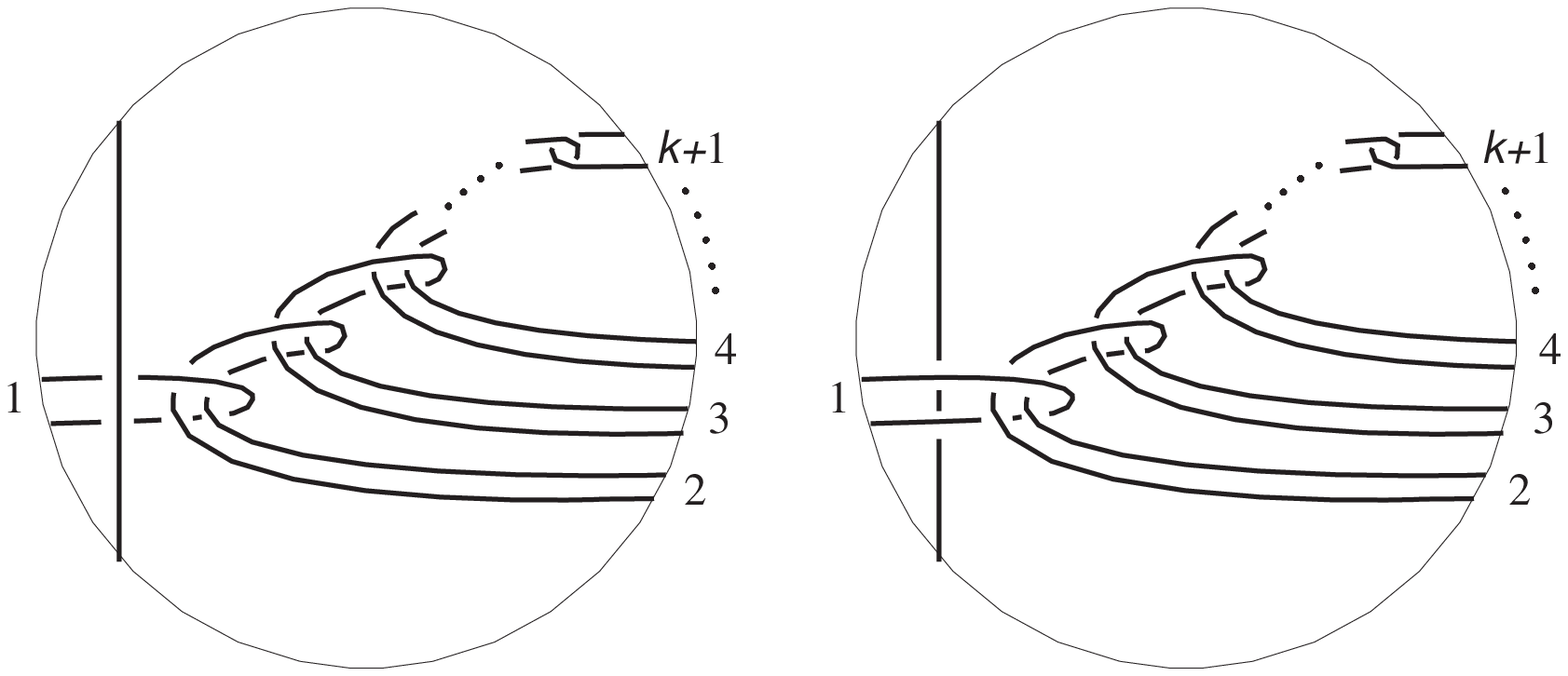} 

Fig. 3 
\end{center}

Let $L_1=K_{11}\cup\cdots\cup K_{1n}$ and 
$L_2=K_{21}\cup\cdots\cup K_{2n}$ be oriented ordered 
$n$-component links (or ordered $n$-component tangle) 
in $S^3$, and 
let $(T_{1},T_{2})$ be a local move. Suppose that 
there is an orientation preserving embedding 
$h:B^{3}\rightarrow S^{3}$ such that 
$(h^{-1}(L_1),h^{-1}(L_2))\cong(T_1,T_2)$ and 
$L_{1}-h(B^{3})=L_{2}-h(B^{3})$ together with orientations 
and orders. 
Then $L_2$ is said to be  {\em obtained from 
$L_1$ by $(T_1,T_2)$}. 
And the set $\{l|h(B^3)\cap K_{1l}\neq\emptyset\}
(=\{l|h(B^3)\cap K_{2l}\neq\emptyset\})$ is the index of 
$(T_{1},T_{2})$. 
We call a (one-branched) $C_{k}$-move ($k<n$) is 
{\em $($one-branched$)$ d-$C_k$-move}, 
if the number of elements of the index is equal to $k+1$. 

Two oriented ordered links $L_{1}$ and
$L_{2}$ are {\em $C_{k}$-equivalent} 
(resp. {\em d-$C_k$-equivalent, one-branched d-$C_k$-equivalent}) 
if $L_{2}$ is obtained from $L_{1}$ by a finite sequence of $C_{k}$-moves 
(resp. d-$C_{k}$-moves, one-branched d-$C_{k}$-moves) and ambient isotopies. 
These relations are equivalence relations. 

A link is called {\em Brunnian} (or {\em almost trivial} \cite{Mil}) 
if the proper sublinks of it are all trivial. The following proposition 
gives a characterization of Brunnian links.

\medskip
{\bf Proposition 1.1.} {\em Let $L$ be an $n$-component link $(n\geq 2)$.  
The following conditions are mutually 
equivalent. \\
{\rm (1)} $L$ is a Brunnian link.\\
{\rm (2)} $L$ is d-$C_{n-1}$-equivalent to a trivial link.\\
{\rm (3)} $L$ is one-branched d-$C_{n-1}$-equivalent to a trivial link. 
}

\medskip
{\bf Remarks.}  
(1) This proposition implies that an $n$-component Brunnian link is 
ambient isotopic to a {\em band sum} of a trivial link and the 
$n$-component {\em Milnor's link}. We will mention that 
more precisely in the next section (Proposition 2.3).

(2) Propositions 1.1 and 2.3 are also shown by Habiro \cite{Habiro3} 
independently.  

\medskip
Two links are {\em link-homotopic} \cite{Mil}, if they are 
transformed into each other by a finite sequence of self 
crossing change and ambient isotopies. J. Milnor gave a 
classification theorem for Brunnian links up to link homotopy 
by using his $\mu$-invariants \cite{Mil}. In this paper, we have 

\medskip
{\bf Theorem 1.2.} {\em 
Two $n$-component Brunnian links $L$ and $L'$ are 
$C_{n}$-equivalent if and only if 
they are link-homotopic. }

\medskip{\bf Remarks.} 
(1) Since for $k<n$, the $C_n$-equivalence implies the $C_k$-equivalence, 
by Proposition 1.1, $n$-component Brunnian links are 
$C_k$-equivalent to a trivial link if $k<n$. 

(2) A $C_2$-move is equal to the {\em delta-move} defined by 
H. Murakami and Y. Nakanishi \cite{M-N}. The $C_2$-equivalence can be 
classified by the linking number \cite{M-N}. A $C_3$-move is 
equal to the {\em clasp-pass move} defined by Habiro \cite{Habiro}. 
Classifications for 3-component links and for
algebraically split links up to $C_3$-move are 
given by K. Taniyama and the second author \cite{T-Y0}.

\bigskip
\noindent
{\bf 2. Band description of a link}

\bigskip
From now on we consider links up to ambient isotopy of 
$S^3$ and tangles up to ambient isotopy of
$B^3$ relative to $\partial B^3$ without explicit mention. 

In order to prove Proposition 1.1 and Theorem 1.2, 
we need the {\em band description} 
of a link defined in \cite{T-Y}. See also \cite{Suz}, \cite{Yam}, 
\cite{Yas}, \cite{T-Y2}.
A {\em $C_{1}$-link model} is a pair $(\alpha,\beta)$ where 
$\alpha$ is a disjoint union of properly embedded arcs in $B^3$ and $\beta$ 
is a disjoint union of arcs on $\partial B^3$ with 
$\partial \alpha=\partial \beta$ as illustrated in Fig. 4.
Suppose that a $C_{k}$-link model $(\alpha, \beta)$ is defined 
where $\alpha$ is a disjoint union of $k+1$ properly embedded arcs in $B^3$ 
and $\beta$ is a disjoint union of $k+1$ arcs on $\partial B^3$ with 
$\partial \alpha=\partial \beta$ such that $\alpha \cup \beta$ is a disjoint 
union of $k+1$ circles. Let $\gamma$ be a component of $\alpha \cup \beta$ and 
$W$ a regular neighbourhood of $\gamma$ in $(B^3-(\alpha \cup
\beta))\cup\gamma$. Let $V$ be an
oriented solid  torus, $D$ a disk in $\partial V$, $\alpha_{0}$ properly
embedded arcs in $V$ 
and $\beta_{0}$ arcs on $D$ as illustrated in Fig. 5. Let $\psi
:V\rightarrow W$ be an orientation
preserving homeomorphism  such that $\psi (D)=W\cap \partial B^{3}$ and 
$\psi (\alpha_{0}\cup \beta_{0})$ bounds disjoint disks in $B^3$. 
Then we call the pair $((\alpha -\gamma )\cup \psi (\alpha_{0}), 
(\beta -\gamma )\cup\psi (\beta_{0}))$ a {\em $C_{k+1}$-link model.} 

\begin{center} 
\begin{tabular}{cc}
\includegraphics[trim=0mm 0mm 0mm 0mm, width=.17\linewidth]
{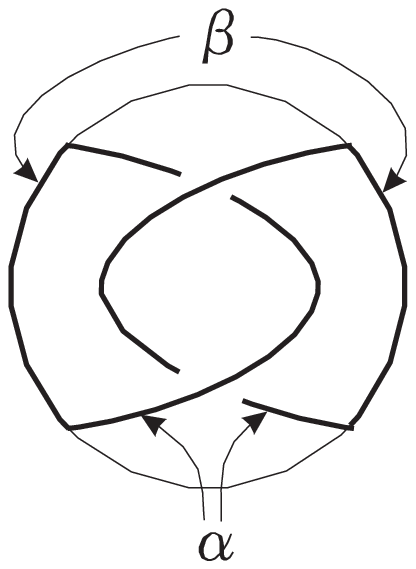} \hspace*{1in} &
\includegraphics[trim=0mm 0mm 0mm 0mm, width=.2\linewidth]
{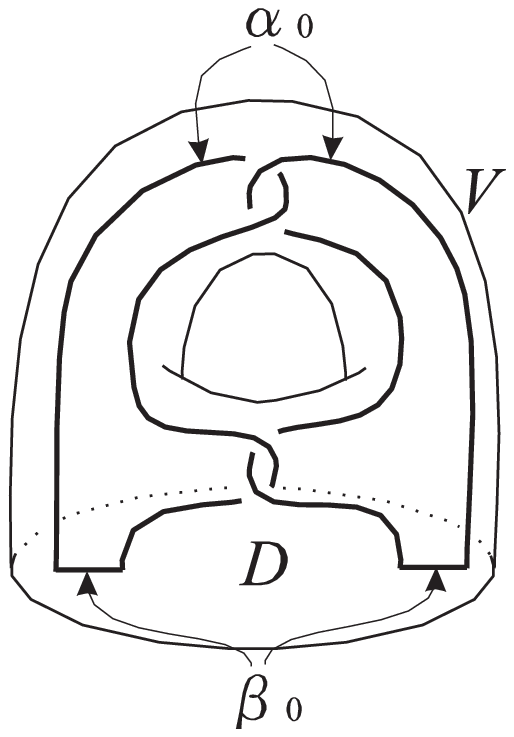}\\
Fig. 4 \hspace*{1in} &
Fig. 5
\end{tabular}
\end{center}

Let $(\alpha_{i},\beta_{i})$ be $C_{\rho(i)}$-link models 
$(i=1,...,l)$. Let $L=K_1\cup \cdots\cup K_n$ 
be an oriented link (resp. a tangle). Let
$\psi_{i}:B^{3}\rightarrow
S^3$ (resp. $\psi_{i}:B^{3}\rightarrow
{\rm int}B^3$)  be an orientation preserving embedding for $i=1,...,
l$ and 
$b_{1,1},b_{1,2},...,b_{1,\rho(1)+1},b_{2,1},b_{2,2},...,b_{2,\rho(2)+1},
...,b_{l,1},b_{l,2},...,b_{l,\rho(l)+1}$
mutually disjoint disks embedded in $S^3$ (resp. $B^3$). Suppose that 
they satisfy the following
conditions;\\
(1) $\psi_{i}(B^{3})\cap \psi_{j}(B^{3})=\emptyset$ if $i\neq j$,\\
(2) $\psi_{i}(B^{3})\cap L=\emptyset$  for each $i$,\\
(3) $b_{i,k}\cap L=\partial b_{i,k}\cap L$ is an arc for each $i,k$,\\
(4) $b_{i,k}\cap (\bigcup_{j=1}^{l} \psi_{j}(B^{3}))=
\partial b_{i,k}\cap \psi_{i}(B^{3})$ is a component of 
$\psi_{i}(\beta_{i})$ for each $i,k$.\\
Let $L'$ be an oriented link (resp. a tangle) defined by
\[
L'=L\cup (\bigcup_{i,k}\partial b_{i,k})\cup 
(\bigcup_{i=1}^{l}\psi_{i}(\alpha_{i})) - 
\bigcup_{i,k}{\rm int}(\partial b_{i,k}\cap L) - 
\bigcup_{i=1}^{l}\psi_{i}({\rm int}\beta_{i}),
\]
where the orientation of $L'$ 
coincides that of $L$ on $L-\bigcup_{i,k}b_{i,k}$ 
if $L$ is oriented.  We call each $b_{i,k}$ a {\it band}. 
We set
${\cal B}_i=((\alpha_i,\beta_i),\psi_i,\{b_{i,1},...,b_{i,\rho(i)+1}\})$
and call ${\cal B}_i$ a {\it
$C_{\rho(i)}$-chord}. 
The set $\{l|K_l\cap\{b_{i,1},b_{i,2},...,b_{i,\rho(i)+1}\}\neq \emptyset\}$ 
is called the {\em index} of ${\cal B}_i$. 
We denote $L'$  by
$L'=\Omega(L;\{{\cal B}_1,...,{\cal B}_l\})$, and say that
$L'$ is a {\em band sum} of
$L$ and chords ${\cal B}_1,...,{\cal B}_l$ or 
a {\em band sum} of
$L$ and $\{{\cal B}_1,...,{\cal B}_l\}$.

A $C_k$-link model $(\alpha,\beta)$ illustrated in 
Fig. 6 is called the {\em one-branched $C_k$-link model}. 
And a $C_k$-chord is called the 
{\em one-branched $C_k$-chord} if its link model is 
the one-branched $C_k$-link model. 
Note that $\alpha\cup\beta \subset S^3$ is the 
$(k+1)$-component Milnor's link \cite[Fig. 7]{Mil}. 

\begin{center}
\includegraphics[trim=0mm 0mm 0mm 0mm, width=.3\linewidth]
{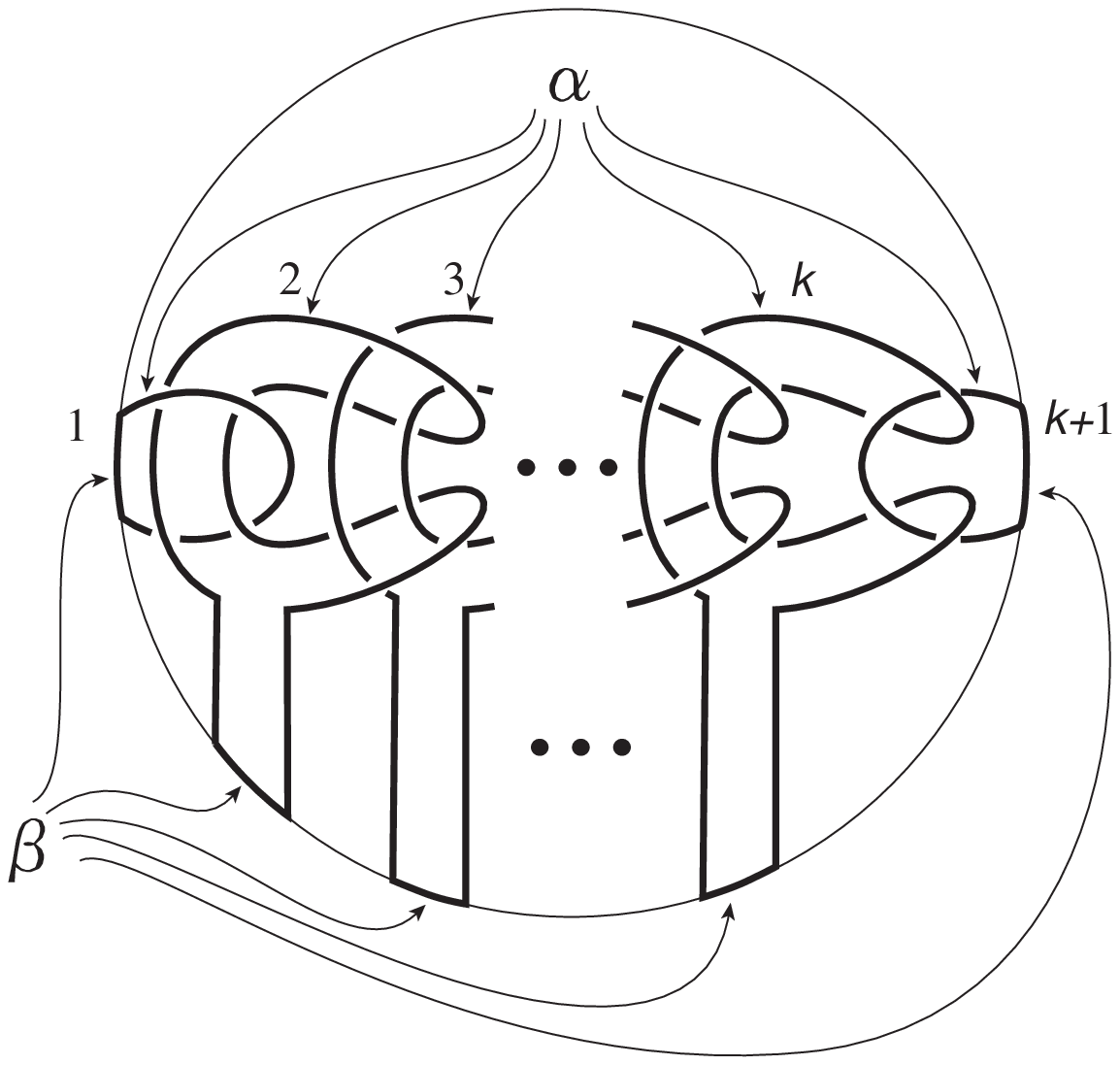} 

Fig. 6
\end{center}

By the arguments similar to that in the proofs of Lemmas 3.6 and 
3.8 in \cite{T-Y}, we have the following two lemmas respectively.

\medskip
\noindent{\bf Lemma 2.1.} (cf. \cite[Lemma 3.6]{T-Y})
{\em A link $L'$ is obtained from $L$ by $($one branched$)$
$C_{k_i}$-moves with indices $I_i\ (i=1,...,m)$ if and only 
if $L'$ is a band sum of $L$ and $($one branched$)$ 
$C_{k_i}$-chords with indices $I_i\ (i=1,...,m)$. $\Box$}

\medskip
\noindent{\bf Lemma 2.2.} (cf. \cite[Lemma 3.8]{T-Y}) 
{\em 
 Let $L=K_1\cup\cdots\cup K_n$, $J=\Omega(L;\{{\cal B}\})$ 
 and $J'=\Omega(L;\{{\cal B}'\})$ be 
 oriented links, where ${\cal B},{\cal B}'$ 
are $C_k$-chords with indices $I$. Suppose that $J$ and $J'$ differ 
locally as illustrated in Fig. $7$ {\rm (a), (b)}, i.e., 
$J'$ is obtained from $J$ by a crossing change 
between $K_i$ and a band of the $C_k$-chord ${\cal B}$. 
Then $J'$ is obtained from $J'$ 
by a $C_{k+1}$-move with index $I\cup\{i\}$. $\Box$ }

\begin{center}
\includegraphics[trim=0mm 0mm 0mm 0mm, width=.4\linewidth]
{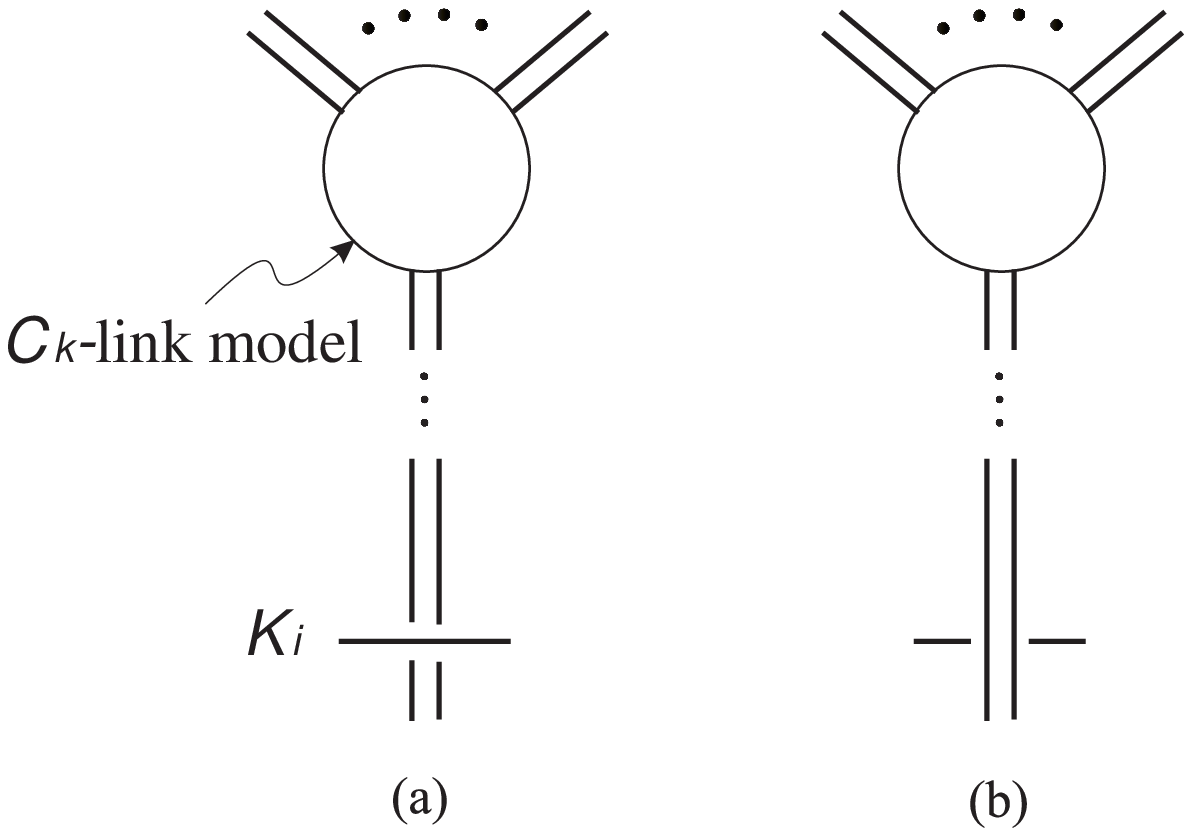}

Fig. 7
\end{center}

Thus by combining Proposition 1.1 and Lemma 2.1, we have 

\medskip
{\bf Proposition 2.3} {\em An $n$-component Brunnian link 
is a band sum of a trivial link and some 
one-branched $C_{n-1}$-chords with indices $\{1,...,n\}$.  
$\Box$}

\bigskip
\noindent
{\bf 3. Proofs of Proposition 1.1 and Theorem 1.2}

\bigskip
A local move $(T_1,T_2)$ is {\em trivial}, if $(T_1,T_2)$ is 
equivalent to a local move $(T_1,T_1)$. 
Let $(T_1,T_2)$ be a local move, and let $t_1,t_2,...,t_k$ and 
$u_1,u_2,...,u_k$ be the components with 
$\partial t_i=\partial u_i\ (i=1,2,...,k)$ of 
$T_1$ and $T_2$ respectively.  
We call $(T_1,T_2)$ a {\em $k$-component Brunnian
local move} ($k\geq 2$), if 
each local move $(T_1-t_i,T_2-u_i)$ is trivial $(i=1,2,...,k)$ 
\cite{T-Y00}. 
It is easy to see that a $C_k$-move is a $(k+1)$-component 
Brunnian local move. 

\medskip
{\bf Proof of Proposition 1.1.} 
 A proof of \lq (2)$\Leftrightarrow$(3)' follows 
 the proof of Lemma 2.2 in \cite{T-Y}. 
Since a $C_{n-1}$-move is an $n$-component Brunnian local move, 
a d-$C_{n-1}$-move preserves link types for any proper sublink 
of an $n$-component link.
Thus we have \lq (2)$\Rightarrow$(1)'.

We will show \lq (1)$\Rightarrow$(2)'. 
Suppose that $L=K_1\cup\cdots\cup K_n$ is an $n$-component Brunnian link. 
Since $L-K_1$ is a trivial link, there is a link diagram 
$\widetilde{L}=\widetilde{K_1}\cup\cdots\cup\widetilde{K_n}$ such 
that $\widetilde{K_2}\cup\cdots\cup\widetilde{K_n}$ has no crossings. 

Since $L-K_2$ is trivial, $L$ is obtained from a trivial link 
by crossing changes between $\widetilde{K_1}$ and $\widetilde{K_2}$. 
Since crossing change between $K_1$ and $K_2$ is a $C_1$-move with index 
$\{1,2\}$, by Lemma 2.1, $L$ is a band sum $\Omega(O;{\bf B}_1)$ 
of a trivial link $O=O_1\cup\cdots\cup O_n$ 
and a set ${\bf B}_1$ of $C_1$-chords  with indices $\{1,2\}$. 

Since $\Omega(O;{\bf B}_1)-O_3(=L-K_3)$ is trivial, 
$L$ is obtained from a trivial link by crossing 
changes between $O_3$ and some bands of the $C_1$-chords. 
By Lemmas 2.1 and 2.2, $L$ is a band sum $\Omega(O;{\bf B}_2)$ 
of a trivial link $O=O_1\cup\cdots\cup O_n$ 
and a set ${\bf B}_2$ of $C_2$-chords  with indices $\{1,2,3\}$. 

Note that $\Omega(O;{\bf B}_2)-O_4$ is trivial. 
Repeating these processes, we have that 
$L$ is a band sum $\Omega(O;{\bf B}_{n-1})$ 
of a trivial link $O$ 
and a set ${\bf B}_{n-1}$ of $C_{n-1}$-chords  with indices 
$\{1,...,n\}$. Hence, by Lemma 2.1, 
$L$ is d-$C_{n-1}$-equivalent a trivial link. $\Box$

\medskip
{\bf Proof of Theorem 1.2.} 
Suppose an $n$-component link $L'$ is obtained from 
an $n$-component link $L$ by a single $C_n$-move. Then 
there is an embedding $h:B^3\rightarrow S^3$ such that 
$(h^{-1}(L),h^{-1}(L'))$ is a $C_n$-move. 
Since the tangles of $C_n$-move are $(n+1)$-component, 
there are two components $s,t$ of $h(B^3)\cap L$ 
(resp. $s',t'$ of $h(B^3)\cap L'$) with $\partial s=\partial s'$ 
$\partial t=\partial t'$ such that $s$ and $t$ are contained in a 
single component of $L$ (resp. $s'$ and $t'$ are contained in a 
single component of $L'$). 
Since $C_{n}$-move is $(n+1)$-component Brunnian local move, 
we may assume $(h^{-1}(L),h^{-1}(L'))$ has a diagram in the unit disk 
such that $h^{-1}(L)$ has no crossings and that 
$h^{-1}(L'-s')$ has no crossings. Since 
$(h^{-1}(L-t),h^{-1}(L'-t'))$ is trivial local move, 
$h^{-1}(L)$ is obtained from $h^{-1}(L')$ 
by crossing changes between $s'$ and $t'$. 
This implies $L$ and $L'$ are link homotopic. 
So the proof of \lq if' part completes.

We will show \lq only if' part. 
Suppose that $L'=K_1\cup\cdots\cup K_n$ is link-homotopic to $L$. 
Since a self crossing change is a $C_1$-move with index $\{i\}$ 
for some $i\in\{1,...,n\}$, then by Lemma 2.1, 
$L'$ is a band sum of $L$ and $C_1$-chords with indices $\{i\}$ 
($i\in\{1,...,n\})$. 
Since $L$ is Brunnian, by Proposition 1.1, $L$ is a band sum 
$\Omega(O;{\bf B}_{n-1})$ of a trivial link $O=O_1\cup\cdots\cup O_n$ 
and a set ${\bf B}_{n-1}$ of $C_{n-1}$-chords with indices 
$\{1,...,n\}$. 
By Sublemma 3.5 in \cite{T-Y}, 
$L'$ is a band sum of $O$ and $C_1$-chords with indices $\{i\}$, 
the set ${\bf B}_{n-1}$ of $C_{n-1}$-chords. 

{\bf Step 1}: Since $L'-K_1$ is a trivial link, by Lemma 2.2, 
$L'$ is obtained from a trivial link by the following moves;\\
(1) $C_1$-moves with indices $\{1\}$ that correspond to $C_1$-chords 
with indices $\{1\}$, \\
(2) $C_2$-moves with indices $\{1,j\}\ (j\geq2)$ 
that correspond to crossing change between $O_1$ and bands of 
$C_1$-chords with indices $\{j\}$,\\
(3) d-$C_{n-1}$-moves correspond to ${\bf B}_{n-1}$. \\
Thus, by Lemma 2.1, we have that $L'$ is a band sum of $O$ and 
a set ${\bf B}_1$ of $C_1$-chords with indices $\{1\}$, 
a set ${\bf B}_2$ of $C_2$-chords with indices $\{1,j\}$, 
and ${\bf B}_{n-1}$. 

{\bf Step 2}
Since $L'-K_2$ is a trivial link, by Lemma 2.2, 
$L'=\Omega(O;{\bf B}_1\cup{\bf B}_2\cup{\bf B}_{n-1})$ is obtained 
from a trivial link by the following moves;\\
(1) $C_2$-moves with indices $\{1,2\}$ that correspond to $C_2$-chords 
with indices $\{1,2\}$ or crossing change between $O_2$ and bands of 
$C_1$-chords with indices $\{1\}$, \\
(2) $C_3$-moves with indices $\{1,2,j\}\ (j\geq3)$ 
that correspond to crossing change between $O_2$ and bands of 
$C_2$-chords with indices $\{1,j\}$, \\
(3) d-$C_{n-1}$-moves correspond to ${\bf B}_{n-1}$. \\
Thus, by Lemma 2.1, we have that $L'$ is a band sum of $O$ and 
a set of $C_2$-chords with indices $\{1,2\}$, 
a set of $C_3$-chords with indices $\{1,2,j\}$, 
and ${\bf B}_{n-1}$. 

Note that $L'-K_3$ is trivial. 
Repeating these processes to Step $n$, we have that $L'$ is obtained from 
a trivial link by $C_n$-moves with indices $\{1,...,n\}$ and 
d-$C_{n-1}$-moves correspond to ${\bf B}_{n-1}$. 
So $L'$ is $C_n$-equivalent to a band sum of $O$ and 
${\bf B}_{n-1}$. 
This band sum may not be same as the band sum $L=\Omega(O;{\bf B}_{n-1})$. 
Since a small regular neighborhood of the bands and the link balls of 
${\bf B}_{n-1}$ are fixed in each step, by Lemma 2.2, these band sums are 
$C_n$-equivalent. This copletes the proof. $\Box$

\bigskip
\noindent
{\bf 4. $C_n$-move and Brunnian local move}

\bigskip
It is known that a $C_n$-move is an $(n+1)$-component Brunnian local move, 
and that a Brunnian local move preserves Vassiliev invariants of order 
$\leq n-1$. Goussarov-Habiro Theorem \cite{Habiro} \cite{Gus} 
implies that two knots are $C_n$-equivalent if and only if their 
values of any Vassiliev invariants of order $\leq n-1$ are same. 
So, if two knots are equivalent up to $(n+1)$-component 
Brunnian local moves, then they are $C_n$-equivalent. 
For links with two or more components, Gussarov-Habiro Theorem 
does not hold. Here, we consider a relation between $C_n$-move and 
$(n+1)$-component Brunnian local moves. 

By the arguments similar to that in the proof of Proposition 1.1, 
we have the following. As we mentioned in the remark of Proposition 1.1, 
this is also shown by Habiro \cite{Habiro3} independently. 

\medskip
{\bf Proposition 4.1.} {\em Let $(T_1,T_2)$ be an 
$n$-component Brunnian local move $(n\geq 2)$. Then $T_2$ is obtained 
from $T_2$ by a finite sequence of $($one-branched$)$ d-$C_{n-1}$-moves. 
} 

\bigskip
\footnotesize{
 }
\end{document}